\documentstyle[twoside,11pt,leqno]{article}

\def\ind{\textrm{ind}} \def\X{{\cal X}}  \def\H{{\cal H}} 
\def\D{{\cal{D}}}
\def\i{\rm{iso}} 
 \def\N{{\texttt{N}}}
\def\C{\texttt{C}} \def\iso{\textrm{iso}}
\def\c{{\cal C}}
\def\h{{Holo_c(\sigma(A))}}

\def\B{B({\cal H})} \def\b{B({\cal X})}
\def\asc{ \textrm{asc}} \def\dsc{ \textrm{dsc}}

\newtheorem{df}{Definition}[section]
\newtheorem{thm}[df]{Theorem} \newtheorem{pro}[df]{Proposition}
\newtheorem{cor}[df]{Corollary} \newtheorem{ex}[df] {Example}
\newtheorem{rema}[df] {Remark} 
\def\sfstp{{\hskip-1em}{\bf.}{\hskip1em}}

\def\subject#1{\renewcommand{\thefootnote}{}\footnote
{AMS(MOS) subject classification (2010). Primary: {#1}}}

\def\keywords#1{\renewcommand{\thefootnote}{}\footnote
{Keywords: {#1}}}

\def\enddemo{\qed \endtrivlist} \expandafter\let\csname
enddemo*\endcsname=\enddemo

\def\qedsymbol{\ifmmode\bgroup\else$\bgroup\aftergroup$\fi
\vcenter{\hrule\hbox{\vrule
height.5em\kern.5em\vrule}\hrule}\egroup}
\def\qed{\ifmmode\else\unskip\nobreak\fi\quad\qedsymbol}

\title{\bf Isolated eigenvalues, poles and compact perturbations of Banach space operators  }
\author{\normalsize B.P. Duggal}
\date{February--18}

\begin{document}

\maketitle \thispagestyle{empty} \vskip-16pt

\subject{47A10, 47A55, 47A53, 47B40} \keywords{Banach space
operator, isolated eigenvalues, left poles, poles, compact perturbations, SVEP,  Fredholm operator, Toeplitz operator, abstract shift condition }

\begin{abstract} Given a  Banach space operator $A$, the isolated eigenvalues $E(A)$ and the poles $\Pi(A)$ (resp., eigenvalues $E^a(A)$ which are isolated points of the approximate point spectrum and the left ploles $\Pi^a(A)$)  of the spectrum of $A$ satisfy $\Pi(A)\subseteq E(A)$ (resp., $\Pi^a(A)\subseteq E^a(A)$), and the reverse inclusion holds if and only if $E(A)$ (resp., $E^a(A)$) has empty intersection with the B-Weyl spectrum (resp., upper B-Weyl spectrum) of $A$. Evidently $\Pi(A)\subseteq E^a(A)$, but no such inclusion exists for $E(A)$ and $\Pi^a(A)$. The study of identities $E(A)=\Pi^a(A)$ and $E^a(A)=\Pi(A)$, and their stability under perturbation by commuting Riesz operators, has been of some interest in the recent past. This paper studies the stability of these identities under perturbation by (non-commuting) compact operators. Examples of analytic Toeplitz operators and operators satisfying the abstract shift condition are considered.
\end{abstract}


\section {\sfstp Introduction} Let $\b$ (resp., $\B$) denote the
algebra of operators, equivalently bounded linear transformations,
on a complex infinite dimensional Banach (resp., Hilbert) space into
itself. For an operator $A\in\b$, let $\sigma(A)$, $\i\sigma(A)$, $\sigma_p(A)$, $\sigma_a(A)$ and $\i\sigma_a(A)$ denote, respectively, the spectrum, the set of
isolated points of $\sigma(A)$, the point spectrum, the approximate point spectrum and the set of
isolated points of $\sigma_a(A)$. Let $\asc(A)$ (resp.,
$\dsc(A)$) denote the ascent (resp., descent) of $A$,
 $A-\lambda$ denote $A-\lambda I$, $\alpha (A-\lambda)=\rm{dim}(A-\lambda)^{-1}(0)$, and let $E(A)$, $E_0(A)$, $E^a(A)$, $E^a_0(A)$, $\Pi(A)$, $\Pi_0(A)$, $\Pi^a(A)$ and $\Pi^a_0(A)$ denote, respectively the sets $E(A)=\{\lambda\in\i\sigma(A):\lambda\in\sigma_p(A)\}$, $E_0(A)=\{\lambda\in E(A): \alpha(A-\lambda)<\infty\}$, $E^a(A)=\{\lambda\in\i\sigma_a(A):\lambda\in\sigma_p(A)\}$, $E_0^a(A)=\{\lambda\in E^a(A): \alpha(A-\lambda)<\infty\}$, $\Pi(A)=\{\lambda\in\sigma(A):\asc(A-\lambda)=\dsc(A-\lambda)<\infty\}$, $\Pi_0(A)=\{\lambda\in\Pi(A):\alpha(A-\lambda)<\infty\}$,
$\Pi^a(A)=\{\lambda\in\sigma_a(A):\asc(A-\lambda)=d <\infty, (A-\lambda)^{d+1}(\X) \hspace{2mm}\mbox{is closed}\}$ and $\Pi_0^a(A)=\{\lambda\in\Pi^a(A):\alpha(A-\lambda)<\infty\}$. The sets $\Pi(A)$, $\Pi^a(A)$, $E(A)$ and $E^a(A)$ 
satisfy the inclusions $\Pi(A)\subseteq\Pi^a(A)\subseteq E^a(A)$ and $\Pi(A)\subseteq E(A)\subseteq E^a(A)$. The reverse inclusions in general do not hold. The reverse inclusions, in particular the properties $$ (P1):\hspace{5mm} E(A)=\Pi^a(A)\hspace{5mm}\mbox{and}\hspace{5mm} (P2):\hspace{5mm} E^a(A)=\Pi(A)$$ and their stability under perturbations by commuting Riesz operators, have been studied in a number of papers in the recent past, amongst them \cite{{AGP},{AAGP},{AGP1}, {D0}, {D1}, {MHR}, {SC}, {HZ}, {ZL}}. It is easily seen that $A\in\b$ satisfies property $(P1)$, $A\in (P1)$, if and only if $E(A)=\Pi^a(A)=\Pi(A)$ and $A\in (P2)$ if and only if $E^a(A)=\Pi(A)=\Pi^a(A)= E(A)$ (thus: $(P2)\Longrightarrow (P1)$). Letting $\sigma_{Bw}(A)$ and $\sigma_{uBw}(A)$ denote, respectively, the B-Weyl and the left (or, upper) B-Weyl spectrum of $A\in\b$, it is seen that $A\in (P1)$ if and only if $E(A)\cap\sigma_{uBw}(A)=\emptyset$ and $A\in (P2)$ if and only if $E^a(A)\cap\sigma_{Bw}(A)=\emptyset$: Left polaroid operators (i.e., operators $A$ for which $\lambda\in\i\sigma_a(A)$ implies $\lambda\in\Pi^a(A)$) satisfy $(P1)$ and a-polaroid operators (i.e., operators $A$ for which $\lambda\in\i\sigma_a(A)$ implies $\lambda\in\Pi(A)$) satisfy $(P2)$ \cite{D1}. The isolated points of (the Weyl spectrum $\sigma_w(A)$ and) the left (or, upper) Weyl spectrum $\sigma_{aw}(A)$ of $A$ play an important role in determining the stability of properties $(P1)$ and $(P2)$ under perturbation by commuting Riesz operators $R\in\b$. Thus, if $\i\sigma_{aw}(A)=\emptyset$, and $\i\sigma_a(A+R)=\i\sigma_a(A)$, then $A\in (Pi)\Longleftrightarrow A+R\in (Pi)$; $i=1,2$  \cite[Theorem 8.5]{D1}.
  
\

This paper considers the preservation of properties $(P1)$ and $(P2)$, and their finite dimensional kernel versions
$$(P1)':\hspace{5mm} E_0(A)=\Pi_0^a(A)\hspace{5mm}\mbox{and}\hspace{5mm} (P2)':\hspace{5mm} E_0^a(A)=\Pi_0(A),$$ under perturbation by (non-commuting) compact operators. We give a number of examples to show that neither of the properties $(P1)$, $(P1)'$, $(P2)$ and $(P2)'$ travels well from $A$ to $A+K$ under perturbation by compact operators $K\in\b$. It is proved that if $\i\sigma_a(A+K)=\i\sigma_a(A)$ and either $\i\sigma_w(A)\cap\{\sigma(A)\setminus\sigma_{Bw}(A)\}=\emptyset$ or $\i\sigma_{aw}(A)\cap\{\sigma_a(A)\setminus\sigma_{uBw}(A)\}=\emptyset$, then  $A\in(P1)\Longrightarrow A\in (P1)'$, and $A\in (P1)\Longrightarrow A+K\in (P1)'$, if and only if $E_0(A+K)\subseteq E_0(A)$;
 $A\in(P2)\Longrightarrow A\in (P2)'$, and $A\in (P2)\Longrightarrow A+K\in (P2)'$, if and only if $E^a_0(A+K)\subseteq E^a_0(A)$. For $A, K$ such that $\i\sigma_a(A+K)=\i\sigma_a(A)$, $\sigma_{Bw}(A)\setminus\sigma_{uBw}(A)=\sigma_{Bw}(A+K)\setminus\sigma_{uBw}(A+K)$ and $\i\sigma_a(A)\cap\{\sigma_{Bw}(A+K)\setminus\sigma_{Bw}(A)\}=\emptyset$, a sufficient condition for $A\in (Pi)$ implies $A+K\in (Pi)$,
$i=1,2$, is that $\i\sigma_a(A)\cap\sigma_{uBw}(A)=\emptyset$. Analytic Toeplitz operators $A_f\in B(H^2(\partial{\D}))$, and operators $A\in\b$ satisfying the abstract shift condition (such that $A$ is non-quasinilpotent and non-invertible), satisfy properties $(P1)$ and $(P2)$. We prove that a sufficient condition for $A_f+K\in (P1)\vee (P2)$ is $E^a(A+K)\cap\sigma_w(A)=\emptyset$, and a necessary and sufficient condition for $A+K\in (P1)\vee (P2)$ is $\sigma_a(A+K)\cap\{\sigma_w(A)\setminus\sigma_{aw}(A)\}=\emptyset$. 

\

The plan of this paper is as follows. After introducing  (most of) our notation and terminology in Section 2, we prove some complementary results on polaroid type operators, and a functional calculus for such operators, in Section 3. Section 4 is devoted to proving our main results, and Section 5 considers  examples of analytic Toeplitz operators and operators which satisfy the abstract shift condition. 

\

\section {\sfstp  Notation and terminology}  In addition to the (explained) notation and terminology already introduced, 
we shall use the following further notation and 
terminology.We shall use $\C$ to denote the complex plane, and $S^{\c}$ to denote the complement of the subset $S$ of $\C$ in $\C$. (Thus $\sigma_w(A)^{\c}=\C\setminus\sigma_w(A)$.) We use $\D(0,r)$ to denote the open disc (in $\C$)) of radius $r$ centered at $0$, $\D$ to denote (the open) unit disc, $\overline{\D}$ to denote the closure of $\D$ and $\partial{\D}$ to denote the boundary of $\D$. An operator $A\in\b$ has SVEP, {\em the single-valued
extension property}, at $\lambda_0\in\C$ if for every open
disc ${D}_{\lambda_0}$ centered at $\lambda_0$ the only
analytic function $f:{D}_{\lambda_0}\longrightarrow \X$
satisfying $(A-\lambda)f(\lambda)=0$ is the function $f\equiv 0$.
Every operator $A$ has SVEP at points in its resolvent set
$\rho(A)=\C\setminus \sigma(A)$ and  on the boundary $\partial\sigma(A)$
of the spectrum $\sigma(A)$. {\em We say that $T$ has SVEP on a set
$S$ if it has SVEP at every $\lambda\in S$.} The {\em ascent of
$A$}, $\asc(A)$ (resp. {\em descent of $A$}, $\dsc(A)$), is the
least non-negative integer $n$ such that $A^{-n}(0)=A^{-(n+1)}(0)$
(resp., $A^n(\X)=A^{n+1}(\X)$): If no such integer exists, then
$\asc(A)$, resp. $\dsc(A)$, $=\infty$. It is well known  that
$\asc(A)<\infty$ implies $A$ has SVEP at $0$, $\dsc(A)<\infty$
implies $A^*$ ($=$ the dual operator) has SVEP at $0$, finite ascent
and descent for an operator implies their equality, and that a point
$\lambda\in\sigma(A)$ is a pole (of the resolvent) of $A$ if and
only if $\asc(A-\lambda)=\dsc(A-\lambda)<\infty$ (see \cite{{A},
{H}, {LN}, {TL}}).

\

The operator  $A\in \b$ is: {\em left semi--Fredholm at $\lambda\in\C$},
denoted $\lambda\in\Phi_{+}(A)$ or $A-\lambda\in\Phi_{+}(\X)$, if
$(A-\lambda)\X$ is closed and the deficiency
 index
 $\alpha(A-\lambda)<\infty$; {\em
right semi--Fredholm at $\lambda\in\C$}, denoted $\lambda\in\Phi_{-}(A)$ or
$A-\lambda\in\Phi_{-}(\X)$, if $\beta(A-\lambda)=\rm{dim}(\X/
(A-\lambda)(\X))<\infty$.
 $A$ is semi--Fredholm, $\lambda\in\Phi_{sf}(A)$ or  $A-\lambda\in\Phi_{sf}(\X)$,
 if $A-\lambda$ is either left or right semi--Fredholm,  and $A$ is Fredholm,
 $\lambda\in\Phi(A)$ or $A-\lambda\in\Phi(\X)$, if $A-\lambda$ is both left and
 right semi--Fredholm. The index of a semi--Fredholm operator  is
the integer, possibly infinite, $\ind(A)=\alpha(A)-\beta(A)$.
Corresponding to these classes of one sided Fredholm operators, we
have the following spectra: The {\em left Fredholm spectrum}
$\sigma_{ae}(A)$ of $A$ defined by
$\sigma_{ae}(A)=\{\lambda\in\sigma(A):A-\lambda\notin\Phi_{+}(\X)\}$,
and the {\em right Fredholm spectrum} $\sigma_{se}(A)$ of $A$
defined by
$\sigma_{se}(A)=\{\lambda\in\sigma(A):A-\lambda\notin\Phi_{-}(\X)\}$.
The {\em Fredholm spectrum} $\sigma_e(A)$ of $A$ is the set
$\sigma_e(A)=\sigma_{ae}(A)\cup\sigma_{se}(A)$. $A\in\b$ is
 Weyl (resp. a-Weyl)  if it is  Fredholm with $\ind(A)=0 $ (resp., 
 if it is left Fredholm with $\ind(A)\leq 0$). It is well known
that {\em a semi- Fredholm operator $A$ (resp., its dual operator
$A^*$) has SVEP at a point $\lambda$ if and only if
$\asc(A-\lambda)<\infty$ (resp., $\dsc(A-\lambda)<\infty$) \cite[Theorems 3.16, 3.17]{A}; furthermore, if $A-\lambda$ is Weyl
 , i.e., if $\lambda\in\Phi(A)$ and
$\ind(A-\lambda)=0$, then $A$ has SVEP at $\lambda$ implies
$\lambda\in\i\sigma(A)$ with
$\asc(A-\lambda)=\dsc(A-\lambda)<\infty$ (resp., if $A-\lambda$ is   a-Weyl
, i.e., if $\lambda\in\Phi_+(A)$ and
$\ind(A-\lambda)\leq 0$, then $A$ has SVEP at $\lambda$ implies
$\lambda\in\i\sigma_a(A)$ with
$\asc(A-\lambda)<\infty$).} The Weyl (resp., the left
or approximate Weyl) spectrum of $A$ is the set \begin{eqnarray*} &
& \sigma_{w}(A)=\{\lambda\in\sigma(A):
\lambda\notin\Phi(A)\hspace{2mm}\mbox{or}\hspace{2mm}\ind(A-\lambda)\neq
0\}\\ & & (\sigma_{aw}(A)=\{\lambda\in\sigma_a(A):
\lambda\notin\Phi_{+}(A)\hspace{2mm}\mbox{or}\hspace{2mm}\ind(A-\lambda)>
0\}).\end{eqnarray*} .

A generalization of  Fredholm  and Weyl spectra is obtained as follows. 
An operator  {\em $A\in\b$ is semi
	B-Fredholm} if there exists an integer $n\geq
1$ such that $A^n(\X)$ is closed and the induced operator
$A_{[n]}=A|_{A^n(\X)}$, $A_{[0]}=A$, is semi Fredholm (in the usual sense). It is seen
that if $A_{[n]}\in\Phi_{\pm}(\X)$ for an integer $n\geq 1$, then
$A_{[m]}\in\Phi_{\pm}(\X)$ for all integers $m\geq n$, and one may
unambiguously define the index of $A$ by
$\ind(A)=\alpha(A)-\beta(A)$ ($=\ind(A_{[n]})$) (see \cite{D00} and \cite{AGP1} for relevant references). Upper (or, left)
semi B-Fredholm, lower (or, right) semi B-Fredholm and B-Fredholm spectra of $A$
are then  the sets 
\vskip2pt $ \sigma_{uBf}(A)  = 
\{\lambda\in\sigma(A): A-\lambda\hspace{2mm}\mbox{is not upper semi
	B-Fredholm}\}$,

\vskip2pt $ \sigma_{lBf}(A)  = 
\{\lambda\in\sigma(A): A-\lambda\hspace{2mm}\mbox{is not lower semi
	B-Fredholm}\}$, and 
\vskip2pt $\sigma_{Be}(A)=\sigma_{uBf}(A)\cup\sigma_{lBf}(A)$.

\vskip2pt

\noindent Letting
\vskip2pt $\sigma_{Bw}(A)=\{\lambda\in\sigma(A):\lambda\in\sigma_{Be}(A)\hspace{2mm}
	\mbox{or}\hspace{2mm}\ind(A-\lambda)\neq
	0\},$
	\vskip2pt $\sigma_{uBw}(A)=\{\lambda\in\sigma_a(A): \lambda
	\in\sigma_{uBf}(A)\hspace{2mm}\mbox{or}\hspace{2mm}
	\ind(A-\lambda)\not\leq 0\},$ 
	\vskip2pt $	\sigma_{lBw}(A)=\{\lambda\in\sigma_s(A): \lambda
	\in\sigma_{lBf}(A)\hspace{2mm}\mbox{or}\hspace{2mm}
	\ind(A-\lambda)\not\geq 0\}$ 
	\vskip4pt\noindent denote, respectively, the {\em the
	B-Weyl, the upper B-Weyl and the lower B-Weyl spectrum of $A$}, we have
$\sigma_{Bw}(A)=\sigma_{uBw}(A)\cup\sigma_{lBw}(A)$, and $\sigma_{uBw}(A)=\sigma_{lBw}(A^*)$. Just as in the case of Weyl and a-Weyl operators, if $A$ has SVEP at $\lambda\in\sigma(A)$ and $A-\lambda$ is B-Weyl, then $\asc(A-\lambda)=\dsc(A-\lambda)<\infty$ and $\lambda\in\i\sigma(A)$ (resp., if $A$ has SVEP at $\lambda\in\sigma_a(A)$ and $A-\lambda$ is upper B-Weyl, then $\asc(A-\lambda)<\infty$ and $\lambda\in\i\sigma_a(A)$) \cite{D00}. 

\vskip4pt

We say in the following that $A\in\b$ is {\em polaroid} (resp., {\em finitely polaroid}) if $\i\sigma(A)\subseteq\Pi(A)$ (resp., $\i\sigma(A)\subseteq \Pi_0(A)$), {\em left polaroid} (resp., {\em finitely left polaroid}) if $\i\sigma_a(A)\subseteq\Pi^a(A)$ (resp., $\i\sigma_a(A)\subseteq \Pi^a_0(A)$), {\em a-polaroid} (resp., {\em finitely a-polaroid}) if $\i\sigma_a(A)\subseteq\Pi(A)$ (resp., $\i\sigma_a(A)\subseteq \Pi_0(A)$). It is clear that  a-polaroid operators are polaroid,  $\Pi_0(A)\subseteq\Pi^a_0(A)\subseteq\Pi^a(A)$ and $\Pi_0(A)\subseteq\Pi(A)\subseteq\Pi^a(A)$ \cite{{AGP1}, {D00}, {LZ}}.

\section {\sfstp  Polaroid operators and compact perturbations} Given operators $A, K\in\b$ with $K$ compact, it is clear that 
\begin{eqnarray*}
 & & A+K\hspace{2mm}\mbox{is finitely polaroid} \Longleftrightarrow \iso\sigma(A+K)\subseteq\Pi_0(A+K)\\ &\Longleftrightarrow& 
\iso\sigma(A+K)\cap\sigma_w(A)=\emptyset;\\
	& & A+K\hspace{2mm}\mbox{is finitely left polaroid}\Longleftrightarrow \iso\sigma_a(A+K)\subseteq\Pi_0^a(A+K)\\ &\Longleftrightarrow& 
	\iso\sigma_a(A+K)\cap\sigma_{aw}(A)=\emptyset,\hspace{2mm}\mbox{and}\\
	& & A+K\hspace{2mm}\mbox{is finitely a-polaroid}\Longleftrightarrow \iso\sigma_a(A+K)\subseteq\Pi_0(A+K)\\ &\Longleftrightarrow& 
	\iso\sigma_a(A+K)\cap\sigma_{w}(A)=\emptyset.\end{eqnarray*}
A version of these observations extends to polaroid, left polaroid and a-polaroid operators. Recall from \cite[Section 3]{D1} that $$\sigma_{Bw}(A)=\sigma_w(A)\setminus{\Phi^{\iso}_{Bw}(A)}\hspace{2mm}\mbox{and}\hspace{2mm} \sigma_{uBw}(A)=\sigma_{aw}(A)\setminus{\Phi^{\iso}_{uBw}(A)},$$ where $$\Phi^{\iso}_{Bw}(A)=\iso\sigma_w(A)\cap\sigma_{Bw}(A)^{\c}
\hspace{2mm}\mbox{and}\hspace{2mm} \Phi^{\iso}_{uBw}(A)=\iso\sigma_{aw}(A)\cap\sigma_{uBw}(A)^{\c}.$$ Recall also that for every $\lambda\in\iso\sigma(A+K)$, either $\lambda\in\sigma_w(A+K)=\sigma_w(A)$ or $\lambda\in\sigma_w(A+K)^{\c}=\sigma_w(A)^{\c}$ (similary, 
for every $\lambda\in\iso\sigma_a(A+K)$, either $\lambda\in\sigma_{aw}(A+K)=\sigma_{aw}(A)$ or $\lambda\in\sigma_{aw}(A+K)^{\c}=\sigma_{aw}(A)^{\c}$). Hence, since $\sigma_w(A+K)^{\c}\cap\sigma_{Bw}(A+K)=\emptyset=\sigma_{aw}(A+K)^{\c}
\cap\sigma_{uBw}(A+K)$, 
 \begin{eqnarray*} \iso\sigma(A+K)\cap\sigma_{Bw}(A+K) &=& \{\iso\sigma(A+K)\cap\sigma_w(A+K)\}\cap\sigma_{Bw}(A+K)\\ &\subseteq& \iso\sigma_w(A+K)\cap\sigma_{Bw}(A+K)\\ &=& \iso\sigma_w(A)\cap\sigma_{Bw}(A+K)\end{eqnarray*} and 
 \begin{eqnarray*} \iso\sigma_a(A+K)\cap\sigma_{uBw}(A+K) &=& \{\iso\sigma_a(A+K)\cap\sigma_{aw}(A+K)\}\cap\sigma_{uBw}(A+K)\\ &\subseteq& \iso\sigma_{aw}(A+K)\cap\sigma_{uBw}(A+K)\\ &=& \iso\sigma_{aw}(A)\cap\sigma_{uBw}(A+K).\end{eqnarray*}

 The following theorem, which gives a necessary and sufficient condition for the perturbation of an operator by a compact operator to be 
 polaroid (left polaroid, a-polaroid), improves \cite[Theorem 6.4]{DK}. Let $[\i\sigma_w(A)]_K$ and $[\i\sigma_{aw}(A)]_K$ denote, respectively, the sets $$ [\i\sigma_w(A)]_K=\{\lambda\in\i\sigma_w(A)=\i\sigma_w(A+K):\lambda\in\i\sigma(A+K)\}$$ and 
 $$[\i\sigma_{aw}(A)]_K=\{\lambda\in\i\sigma_{aw}(A)=\i\sigma_{aw}(A+K):\lambda\in\i\sigma_a(A+K)\}.$$ Clearly, $\i\sigma(A+K)\cap\sigma_{Bw}(A+K)=\emptyset\Longleftrightarrow [\i\sigma_w(A)]_K\cap\sigma_{Bw}(A+K)=\emptyset$ and $\i\sigma_a(A+K)\cap\sigma_{uBw}(A+K)=\emptyset\Longleftrightarrow [\i\sigma_{aw}(A)]_K\cap\sigma_{uBw}(A+K)=\emptyset$.
\begin{thm}
	\label{thm1} If $A,K\in\b$, then:

\noindent (i) \begin{eqnarray*}A+K\hspace{2mm}\mbox{is polaroid} &\Longleftrightarrow& \iso\sigma(A+K)\cap\sigma_{Bw}(A+K)=\emptyset\\ &\Longleftrightarrow& [\iso\sigma_w(A)]_K\cap\sigma_{Bw}(A+K)=\emptyset.\end{eqnarray*}	
(ii) \begin{eqnarray*}A+K\hspace{2mm}\mbox{is left polaroid} &\Longleftrightarrow& \iso\sigma_a(A+K)\cap\sigma_{uBw}(A+K)=\emptyset\\ &\Longleftrightarrow& [\iso\sigma_{aw}(A)]_K\cap\sigma_{uBw}(A+K)=\emptyset.		
	\end{eqnarray*}
(iii) \begin{eqnarray*}A+K\hspace{2mm}\mbox{is a-polaroid} &\Longleftrightarrow& \iso\sigma_a(A+K)\cap\sigma_{Bw}(A+K)=\emptyset.
	\end{eqnarray*}
\end{thm} 	
\begin{demo} Start by observing that $$\Pi(A+K)=\Pi_0(A+K)\cup\Pi_{\infty}(A+K)$$ and $$\Pi^a(A+K)=\Pi^a_0(A+K)\cup\Pi^a_{\infty}(A+K),$$ 
	where $$\Pi_{\infty}(A+K)=\iso\sigma(A+K)\cap\{\sigma_w(A+K)\setminus\sigma_{Bw}(A+K)\}$$ and 
	$$\Pi_{\infty}^a(A+K)=\iso\sigma_a(A+K)\cap\{\sigma_{aw}(A+K)\setminus\sigma_{uBw}(A+K)\}.$$
	
\noindent (i) We have:\begin{eqnarray*} \iso\sigma(A+K) &=& \{\iso\sigma(A+K)\cap\sigma_w(A+K)^{\c}\}\cup\{
		\i\sigma(A+K)\\ & & \cap(\sigma_w(A+K)\setminus\sigma_{Bw}(A+K))\}\cup\{\iso\sigma(A+K)\cap\sigma_{Bw}(A+K)\}\\ 
		&=& \Pi_0(A+K)\cup\Pi_{\infty}(A+K)\cup\{\iso\sigma(A+K)\cap\sigma_{Bw}(A+K)\}\\ &=& 
		\Pi(A+K)\cup\{\iso\sigma(A+K)\cap\sigma_{Bw}(A+K)\}\\ &=& \Pi(A+K)\cup\{[\iso\sigma_w(A)]_K\cap\sigma_{Bw}(A+K)\},\end{eqnarray*} which implies 
	$$\iso\sigma(A+K)=\Pi(A+K)\Longleftrightarrow \iso\sigma(A+K)\cap\sigma_{Bw}(A+K)=\emptyset\Longleftrightarrow [\i\sigma_w(A)]_K\cap\sigma_{Bw}(A+K)=\emptyset.$$	
\noindent (ii) Again: \begin{eqnarray*} \iso\sigma_a(A+K) & = & \{\iso\sigma_a(A+K)\cap\sigma_{aw}(A+K)^{\c}\}\cup\{
		\iso\sigma_a(A+K)\\ & & \cap(\sigma_{aw}(A+K)\setminus\sigma_{uBw}(A+K))\} \cup\{\i\sigma_a(A+K)\cap\sigma_{uBw}(A+K)\}\\ 
		&=& \Pi_0^a(A+K)\cup\Pi_{\infty}^a(A+K)\cup\{\iso\sigma_a(A+K)\cap\sigma_{uBw}(A+K)\}\\ &=& 
		\Pi^a(A+K)\cup\{\iso\sigma_a(A+K)\cap\sigma_{uBw}(A+K)\}\\ &=& \Pi^a(A+K)\cup\{[\i\sigma_{aw}(A)]_K\cap\sigma_{uBw}(A+K)\},	\end{eqnarray*} which implies 
	$$\iso\sigma_a(A+K)=\Pi^a(A+K)\Longleftrightarrow \iso\sigma_a(A+K)\cap\sigma_{uBw}(A+K)=\emptyset\Longleftrightarrow[\iso\sigma_{aw}(A)]_K\cap\sigma_{uBw}(A+K)=\emptyset.$$	
	\noindent (iii) Finally \begin{eqnarray*} \iso\sigma_a(A+K) & = & \{\iso\sigma_a(A+K)\cap\sigma_w(A+K)^{\c}\}\cup\{
		\iso\sigma_a(A+K)\\ & &\cap(\sigma_w(A+K)\setminus\sigma_{Bw}(A+K))\} \cup\{\iso\sigma_a(A+K)\cap\sigma_{Bw}(A+K)\}\\ 
		&=& \Pi_0(A+K)\cup\Pi_{\infty}(A+K)\cup\{\iso\sigma_a(A+K)\cap\sigma_{Bw}(A+K)\}\\ &=& 
		\Pi(A+K)\cup\{\iso\sigma_a(A+K)\cap\sigma_{Bw}(A+K)\},	\end{eqnarray*} which implies 
	$$\iso\sigma_a(A+K)=\Pi(A+K)\Longleftrightarrow \iso\sigma_a(A+K)\cap\sigma_{Bw}(A+K)=\emptyset.$$
This completes the proof. \end{demo}
	
	\begin{rema}
		\label{rema0} {\em Commuting Riesz operators. Translated to operators $A\in\b$ and Riesz operators $R\in\b$ such that 
			$[A,R]=AR-RA=0$ and $\sigma_{Bw}(A+R)=\sigma_{Bw}(A)$ (resp., $\sigma_{uBw}(A+R)=\sigma_{uBw}(A)$) Theorem \ref{thm1} 
			implies that: {\em $A+R$ is polaroid if and only if $\iso\sigma(A)\cap\sigma_{Bw}(A)=\emptyset$, equivalently if and only if $A$ is polaroid (resp., $A+R$ is left polaroid if and only if $\iso\sigma_{a}(A)\cap\sigma_{uBw}(A)=\emptyset$, equivalently if and only if $A$ is left polaroid).} An important example of a class of operators satisfying the above spectral hypotheses is that of operators $F\in\b$ satisfying $[A,F]=0$ and $F^n$ is finite rank for some natural number $n$ \cite[Proposition 3.3]{D1}.}
	\end{rema}

 \noindent {\bf Functional Calculus} Given $A\in\b$, let ${Holo}(\sigma(A))$ denote the set of functions $f$ which are holomorphic in a neighbourhood of $\sigma(A)$, and let $\h$ denote those $f\in{Holo}(\sigma(A))$ which are non-constant 
 on the connected components of $\sigma(A)$. If we let $\sigma_D(A)$ denote the Drazin spectrum of $A$, $$\sigma_D(A)=\{\lambda\in\sigma(A):\asc(A-\lambda)\neq\dsc(A-\lambda)\},$$ then $\sigma_D(A)$ satisfies the spectral mapping theorem $$\sigma_D(f(A))=f(\sigma_D(A)), f\in\h;$$ the left Drazin spectrum $\sigma_{lD}(A)$ of $A$, 
 \begin{eqnarray*} \sigma_{lD}(A) &=& \{\lambda\in\sigma_a(A):\mbox{there does not exist an integer}\hspace{2mm}p\geq 1\\ & & 
 	\hspace{2mm}\mbox{such that}\hspace{2mm}\asc(A-\lambda)\leq p\hspace{2mm}\mbox{and}\hspace{2mm}(A-\lambda)^{p+1}(\X)\hspace{2mm}\mbox{is closed}\}, \end{eqnarray*} 
 also satisfies a similar spectral mapping theorem: $$ \sigma_{lD}(f(A)=f(\sigma_{lD}(A)), f\in\h$$ \cite{Mu}. 
 It is straightforward to see that:
 \vskip6pt  $A$ is polaroid if and only if $\iso\sigma(A)\cap\sigma_D(A)=\emptyset$;
 \vskip3pt   $A$ is left polaroid if and only if $\iso\sigma_a(A)\cap\sigma_{lD}(A)=\emptyset$, and
 \vskip3pt   $A$ is a-polaroid if and only if $\iso\sigma_a(A)\cap\sigma_D(A)=\emptyset$.
 \vskip6pt\noindent It is well known (see, for example, \cite[Lemma 4.1]{DDHZ}) that if $f\in\h$, then $\iso\sigma(f(A))=f(\iso\sigma(A))$. Hence, for $f\in\h$, \begin{eqnarray*} f(A)\hspace{2mm}\mbox{is polaroid} 
 	&\Longleftrightarrow& \iso\sigma(f(A))\cap\sigma_D(f(A))=\emptyset\\ &\Longleftrightarrow& f(\iso\sigma(A))\cap f(\sigma_D(A))=
 	f(\iso\sigma(A)\cap\sigma_D(A))=\emptyset\\ &\Longleftrightarrow& A \hspace{2mm}\mbox{is polaroid}.\end{eqnarray*}
 (See \cite{DDHZ} for other alternative arguments.) This argument does not extend to left polaroid operators (for the reason that the spectral mapping theorem fails for $\iso\sigma_a(A)$). However, given a $\lambda\in\iso\sigma_a(f(A))$ for an $f\in\h$, there always exists a $\mu\in\iso\sigma_a(A)$ such that $f(\mu)=\lambda$. Hence 
 \begin{eqnarray*} A\hspace{2mm}\mbox{is left polaroid} 
 	&\Longrightarrow& f(\iso\sigma_a(A)\cap\sigma_{lD}(A))=\emptyset\\ &\Longleftrightarrow& \{f(\mu):\mu\in(\iso\sigma_a(A)\}\cap f(\sigma_{lD}(A))=\emptyset\\ &\Longleftrightarrow& \{\lambda\in\iso\sigma_a(f(A):\lambda=f(\mu),\mu\in\iso\sigma_a(A)\}\cap\sigma_{lD}(f(A))=\emptyset\\
 	&\Longrightarrow& \iso\sigma_a(f(A))\cap\sigma_{lD}(f(A))=\emptyset\\ &\Longleftrightarrow& f(A) \hspace{2mm}\mbox{is left polaroid}.\end{eqnarray*} For the reverse implication, a hypothesis guaranteeing $ f(\iso\sigma_a(A))=iso\sigma_a(f(A))$, such as 
 $f$ is injective or $iso\sigma_a(A)\subseteq\iso\sigma(A)$, is required. It is clear that $A$ is a-polaroid implies $\iso\sigma_a(A)=\iso\sigma(A)$. Hence $$ A\hspace{2mm}\mbox{is a-polaroid}\Longrightarrow f(A)\hspace{2mm}\mbox{is polaroid}\Longrightarrow A\hspace{2mm}\mbox{is polaroid}.$$ Combining with Theorem \ref{thm1}, we have:
 
 \begin{cor}
 	\label{cor0} Given operators $A,K\in\b$ with $K$ compact, and an $f\in\h$:
 	\vskip4pt\noindent (i) $f(A+K)$ is polaroid if and only if $[\iso\sigma_w(A)]_K\cap\sigma_D(A+K)=\emptyset$.
 	\vskip3pt\noindent (ii) If $f$ is injective, then $f(A+K)$ is left polaroid if and only if $[\i\sigma_{aw}(A)]_K\cap\sigma_D(A+K)=\emptyset$.
 	\vskip3pt\noindent (iii) If $f$ is injective, then $f(A+K)$ is a-polaroid if and only if $\iso\sigma_a(A+K)\cap\sigma_D(A+K)=\emptyset$.
 \end{cor}
 \begin{demo} The proof is immediate from Theorem \ref{thm1} once one observes that if an operator $T$ has SVEP at a point $\lambda$, then $\lambda\in\sigma_{Bw}(T)$ (resp., $\lambda\in\sigma_{uBw}(T)$) if and only if $\lambda\in\sigma_D(T)$ (resp., $\lambda\in\sigma_{lD}(T)$).\end{demo}

\section {\sfstp Properties $(P1)$, $(P2)$ and Compact Perturbations} Neither of the properties $(P1)$ and $(P2)$, or their finite kernel versions $$ (P1)'\hspace{1cm} E_0(A)=\Pi^a_0(A)\hspace{3mm}\mbox{and}\hspace{3mm} (P2)'\hspace{1cm} E^a_0(A)=\Pi_0(A),$$ travels well from $A\in\b$ to its perturbation by a compact operator $K\in\b$.
\vskip4pt\noindent \begin{ex}
	\label{ex1} If we let $A=U\oplus Q\in B(\H\oplus\H)$, where $U$ is the forward unilateral shift and $Q$ is an injective 
	compact quasinilpotent operator, then \begin{eqnarray*} & & \sigma_w(A)=\sigma_{Bw}(A)=\overline{\D}, \sigma_{aw}(A)=\partial{\D}\cup\{0\}=\sigma_{uBw}(A), \iso\sigma_w(A)=\emptyset,\\ & & \iso\sigma_{aw}(A)=\{0\}\hspace{2mm}\mbox{and}\hspace{2mm} E(A)=\Pi^a(A)=\emptyset=E^a(A)=\Pi(A).\end{eqnarray*}Let $K\in B(H\oplus\H)$ be the compact operator $K=0\oplus{-Q}$. Then the perturbed operator $A+K=A\oplus 0$ satisfies \begin{eqnarray*} & & \sigma_w(A+K)=\sigma_{Bw}(A+K)=\overline{\D}, \sigma_{aw}(A+K)=\partial{\D}\cup\{0\},\\ & &\sigma_{uBw}(A+K)=\partial{\D}\hspace{2mm}\mbox{and}\hspace{2mm} \iso\sigma_a(A+K)=\iso\sigma_a(A);\end{eqnarray*} hence $$E(A+K)=\emptyset\neq\Pi^a(A+K), E^a(A+K)=\{0\}\neq\Pi(A+K)=\emptyset.$$
	\end{ex}
\begin{ex}
	\label{ex2} Let $A=U\oplus I\in B(\ell^2\oplus\ell^2)$ and $K=0\oplus F\in B(\ell^2\oplus\ell^2)$, where $U\in B(\ell^2)$ is the forward unilateral shift and $F$ is the compact operator $F(x_1,x_2,x_3,\cdots)=(-\frac{x_1}{2}, 0, 0,\cdots)$. Then $$\iso\sigma_w(A)=\iso\sigma_{aw}(A)=\emptyset, \iso\sigma_a(A)=\emptyset\neq\{\frac{1}{2}\}=\iso\sigma_a(A+K)$$ and \begin{eqnarray*} & & E_0(A)=\Pi^a_0(A)=\emptyset=E^a_0(A)=\Pi_0(A),\\ & & \Pi_0(A+K)=\Pi(A+K)=E(A+K)=E_0(A+K)=\emptyset,\\ & & \Pi^a_0(A+K)=\Pi^a(A+K)=E^a(A+K)=E^a_0(A+K)=\{\frac{1}{2}\}.\end{eqnarray*}
\end{ex}
\begin{ex}
	\label{ex3} If we let $A=U\oplus 0\in B(\ell^2\oplus\ell^2)$, where (as before) $U$ is the forward unilateral shift, then \begin{eqnarray*} & &	
	\sigma_w(A)=\sigma_{Bw}(A)=\overline{\D}, \sigma_{aw}(A)=\partial{\D}\cup\{0\}\neq\sigma_{uBw}(A)=\partial{\D},\\ & & \iso\sigma_w(A)=\emptyset\neq\iso\sigma_{aw}(A)=\{0\}\end{eqnarray*} and $$ E_0(A)=\Pi^a_0(A)=\emptyset=E^a_0(A)=\Pi_0(A).$$ 
	Let $K\in B(\ell^2\oplus\ell^2)$ be the compact operator $K=0\oplus Q$, where $Q$ is the compact operator $Q(x_1,x_2,x_3, ...)=(0,\frac{x_2}{2},\frac{x_3}{3},\cdots)$. Then \begin{eqnarray*} & & \sigma_w(A+K)=\sigma_{Bw}(A+K)=\overline{\D}, \sigma_{aw}(A+K)=\partial{\D}\cup\{0\}=\sigma_{uBw}(A+K),\\ & & E(A+K)=E_0(A+K)=\Pi_0(A+K)=\Pi(A+K)=\emptyset, \Pi^a_0(A+K)\\ & & =\Pi^a(A+K)=\{\frac{1}{2},\frac{1}{3}, ...\}, E^a_0(A+K)=E(A+K)=\{0,\frac{1}{2},\frac{1}{3},\cdots\}.	\end{eqnarray*} Evidently, \begin{eqnarray*} & & E_0(A+K)\neq\Pi^a_0(A+K), E(A+K)\neq\Pi^a(A+K),\\ & & E^a_0(A+K)\neq\Pi_0(A+K)\hspace{2mm}\mbox{and}\hspace{2mm} E^a(A+K)\neq\Pi(A+K). \end{eqnarray*}
\end{ex}
The above examples show that neither of the hypotheses $\iso\sigma_a(A)=\iso\sigma_a(A+K)$, $\sigma_{Bw}(A)=\sigma_{Bw}(A+K)$, 
$\sigma_{uBw}(A)=\sigma_{uBw}(A+K)$, $\iso\sigma_{w}(A)=\emptyset$, $\iso\sigma_{aw}(A)=\emptyset$ and (even) $[A,K]=0$ is sufficient to guarantee the transfer of either of the properties $(P1)$, $(P1)'$, $(P2)$ and  $(P2)'$ from $A$ to $A+K$. Recalling, \cite{D1}, $\sigma_{Bw}(A)=\sigma_w(A)\setminus{\Phi^{\iso}_{Bw}(A)}$ and $\sigma_{uBw}(A)=\sigma_{aw}(A)\setminus\Phi^{\iso}_{uBw}(A)$, where ${\Phi^{\iso}_{Bw}(A)}=\iso\sigma_w(A)\cap\sigma_{Bw}(A)^{\c}$ and ${\Phi^{\iso}_{uBw}(A)}=\iso\sigma_{aw}(A)\cap\sigma_{uBw}(A)^{\c}$, we have 
$${\Phi^{\iso}_{Bw}(A)}={\Phi^{\iso}_{Bw}(A+K)}\Longrightarrow \sigma_{Bw}(A)=\sigma_{Bw}(A+K),\hspace{2mm}\mbox{and}$$
$${\Phi^{\iso}_{uBw}(A)}={\Phi^{\iso}_{uBw}(A+K)}\Longrightarrow \sigma_{uBw}(A)=\sigma_{uBw}(A+K).$$ Furthermore, if also $\iso\sigma_a(A)=\iso\sigma_a(A+K)$, then \begin{eqnarray*} \Pi(A) &=& \iso\sigma(A)\cap\sigma_{Bw}(A)^{\c}=\iso\sigma_a(A)\cap\sigma_{Bw}(A)^{\c}\\ &=& \sigma_a(A+K)\cap\sigma_{Bw}(A+K)^{\c}=\Pi(A+K),\end{eqnarray*} and $$\Pi^a(A)=\iso\sigma_a(A)\cap\sigma_{uBw}(A)^{\c}=\iso\sigma_a(A+K)\cap\sigma_{uBw}(A+K)^{\c}=\Pi^a(A+K).$$
This, however, is not enough to warranty
the passage of properties $(P1)$ and $(P2)$ from $A$ to $A+K$.
\begin{ex}
	\label{ex4} Choose $A=Q_1\oplus Q_2\in B(\H\oplus\H)$, where $Q_1$ is an injective compact quasinilpotent operator and $Q_2$ is an injective quasinilpotent such that $Q_2^n(\H)$ is non-closed for all natural numbers $n$. Then
	 \begin{eqnarray*} & & \sigma_w(A)=\sigma_{Bw}(A)=\sigma_{uBw}(A)=\sigma_{aw}(A)=\{0\}, \iso\sigma_w(A)=\iso\sigma_{aw}(A)=\emptyset (\Longrightarrow\\ & & \Phi^{\iso}_{Bw}(A)=\Phi^{\iso}_{uBw}(A)=\emptyset), E(A)=E^a(A)=\Pi^a(A)=\Pi(A)=\emptyset.	 
	 \end{eqnarray*} Now let $K\in B(\H\oplus\H)$ be the compact operator $K=-Q_1\oplus 0$. Then $A+K=0\oplus Q_2$ satisfies 
	  \begin{eqnarray*} & & \iso\sigma_a(A+K)=\iso\sigma_a(A),  \sigma_w(A+K)=\sigma_{Bw}(A+K)=\sigma_{uBw}(A+K)=\\  & & =\sigma_{aw}(A+K)=\{0\}, \iso\sigma_w(A)= \iso\sigma_{aw}(A)=\{0\} (\Longrightarrow \Phi^{\iso}_{Bw}(A+K)=\Phi^{\iso}_{uBw}(A+K)\\ & & =\emptyset), E(A+K)=E^a(A+K)=\{0\}\neq\Pi^a(A+K)=\Pi(A+K)=\emptyset.	 
	  \end{eqnarray*}
\end{ex}
\begin{rema}
	\label{rema2} {\em We note for future reference that the hypothesis $\Phi^{\iso}_{Bw}(A)=\emptyset$ implies $\sigma_{Bw}(A)=\sigma_w(A)$ and the hypothesis $\Phi^{\iso}_{uBw}(A)=\emptyset$ implies $\sigma_{uBw}(A)=\sigma_{aw}(A)$. Furthermore, the hypothesis $\Phi^{\iso}_{uBw}(A)=\emptyset$ also implies $\sigma_{Bw}(A)=\sigma_w(A)$, as the following argument proves. Evidently $$\Phi^{\iso}_{uBw}(A)=\emptyset\Longrightarrow\iso\sigma_{aw}(A)\cap\sigma_{uBw}(A)^{\c}=\emptyset\Longrightarrow 
\iso\sigma_{aw}(A)\subseteq\sigma_{uBw}(A)\subseteq\sigma_{Bw}(A).$$ Take a $\lambda\notin\sigma_{Bw}(A)$. Then $$\lambda\in\sigma_{Bw}(A)^{\c}\subseteq\sigma_{uBw}(A)^{\c}=\sigma_{aw}(A)^{\c}\Longrightarrow\lambda\in\sigma_w(A)^{\c}$$
(since $\lambda\in\sigma_{Bw}(A)^{\c}$ implies $\ind(A-\lambda)=0$, hence $\lambda\notin\sigma_{aw}(A)$ implies $(A-\lambda)(\X)$ is closed and $\ind(A-\lambda)=0$). Thus $\lambda\notin\sigma_w(A)$, consequently $\sigma_w(A)\subseteq\sigma_{Bw}(A)$ (implies $\sigma_w(A)=\sigma_{Bw}(A)$).}
\end{rema}
\begin{thm}
	\label{thm2} Given operators $A, K\in\b$, $K$ compact, such that $\iso\sigma_a(A)=\iso\sigma_a(A+K)$, if either 
	$\Phi^{\iso}_{Bw}(A)=\emptyset=\Phi^{\iso}_{Bw}(A+K)$, or, $\Phi^{\iso}_{uBw}(A)=\emptyset=\Phi^{\iso}_{uBw}(A+K)$, then:
	\vskip4pt\noindent (i) $A\in(P1)\Longrightarrow A\in (P1)'$ and $A\in (P1)\Longrightarrow A+K\in (P1)'$ if and only if $E_0(A+K)\subseteq E_0(A)$.
	\vskip3pt\noindent (ii) $A\in(P2)\Longrightarrow A\in (P2)'$ and $A\in (P2)\Longrightarrow A+K\in (P2)'$ if and only if $E^a_0(A+K)\subseteq E^a_0(A)$.
\end{thm}
\begin{demo} The hypothesis $\Phi^{\iso}_{Bw}(A)=\emptyset=\Phi^{\iso}_{Bw}(A+K)$ implies $$\sigma_w(A)=\sigma_{Bw}(A)=\sigma_{Bw}(A+K)=\sigma_w(A+K)$$ and hence $$\Pi(A)=\iso\sigma(A)\cap\sigma_{Bw}(A)^{\c}=\iso\sigma(A)\cap\sigma_w(A)^{\c}=\Pi_0(A);$$ furthermore, if also  $\iso\sigma_a(A)=\iso\sigma_a(A+K)$, then \begin{eqnarray*} \Pi(A+K) &=& \iso\sigma(A+K)\cap\sigma_{Bw}(A+K)^{\c}\\ &= & \iso\sigma(A+K)\cap\sigma_w(A+K)^{\c}\\ &=&\Pi_0(A+K)= \iso\sigma_a(A+K)\cap\sigma_{Bw}(A+K)^{\c}\\ &=&\iso\sigma_a(A)\cap\sigma_{Bw}(A)^{\c}=\Pi(A)\\ &=&\iso\sigma_a(A)\cap\sigma_w(A)^{\c}=\Pi_0(A).\end{eqnarray*} (Thus $\Pi(A+K)=\Pi_0(A+K)=\Pi_0(A)=\Pi(A)$.)
	Similarly, if  $\Phi^{\iso}_{uBw}(A)=\emptyset=\Phi^{\iso}_{uBw}(A+K)$ and  $\iso\sigma_a(A)=\iso\sigma_a(A+K)$, then 
$$\sigma_{aw}(A)=\sigma_{uBw}(A)=\sigma_{uBw}(A+K)=\sigma_{aw}(A+K)$$ and  
\begin{eqnarray*} \Pi^a(A+K) &=& \iso\sigma_a(A+K)\cap\sigma_{uBw}(A+K)^{\c}\\ &=&\iso\sigma_a(A+K)\cap\sigma_{aw}(A+K)^{\c}= \Pi_0^a(A+K)\\ &=& \iso\sigma_a(A)\cap\sigma_{uBw}(A)^{\c}=\Pi^a(A)\\ &=&\iso\sigma_a(A)\cap\sigma_{aw}(A)^{\c}=\Pi_0^a(A).\end{eqnarray*} (Thus $\Pi^a(A+K)=\Pi^a_0(A+K)=\Pi^a_0(A)=\Pi^a(A)$.)

\

\noindent (i) If $\Phi^{\iso}_{Bw}(A)=\emptyset=\Phi^{\iso}_{Bw}(A+K)$, then \begin{eqnarray*} A\in (P1) &\Longleftrightarrow& E(A)=\Pi^a(A)=\Pi(A)\\ &\Longrightarrow& E_0(A)=\Pi^a_0(A)=\Pi_0(A) (\Longleftrightarrow A\in (P1)')\\ &\Longleftrightarrow& E_0(A)=\Pi^a_0(A)=\Pi_0(A)=\Pi_0(A+K),\end{eqnarray*} and this since \begin{eqnarray*} \Pi^a_0(A+K) &=& \iso\sigma_a(A+K)\cap\sigma_{aw}(A+K)^{\c}\\ &=&\iso\sigma_a(A)\cap\sigma_{aw}(A)^{\c}= \Pi^a_0(A)= \Pi_0(A)=\Pi_0(A+K)\end{eqnarray*} implies $$A\in(P1)\Longrightarrow A\in (P1)'\Longrightarrow E_0(A)=\Pi^a_0(A+K)=\Pi_0(A+K)\subseteq E_0(A+K).$$
Again, if $\Phi^{\iso}_{uBw}(A)=\emptyset=\Phi^{\iso}_{uBw}(A+K)$, then ($\sigma_w(A)=\sigma_{Bw}(A)=\sigma_{Bw}(A+K)=\sigma_w(A+K)$, and) \begin{eqnarray*} A\in (P1) &\Longleftrightarrow& E(A)= (\Pi^a(A) =\Pi(A)=) \Pi^a_0(A)\\ &\Longrightarrow& E_0(A)=\Pi^a_0(A) (\Longleftrightarrow A\in (P1)')\\ &\Longleftrightarrow& E_0(A)=\Pi_0(A)=\Pi^a_0(A)=\Pi^a_0(A+K),\end{eqnarray*} and this since 
$$\Pi^a_0(A+K)=\Pi_0(A)=\iso\sigma_a(A)\cap\sigma_w(A)^{\c}=\iso\sigma_a(A+K)\cap\sigma_{w}(A+K)^{\c}=\Pi_0(A+K)$$ implies 
$$A\in(P1)\Longrightarrow A\in (P1)'\Longrightarrow E_0(A)=\Pi^a_0(A+K)=\Pi_0(A+K)\subseteq E_0(A+K).$$ Hence, in either case,  $A\in(P1)\Longrightarrow A\in (P1)'$ and $A\in (P1)\Longrightarrow A+K\in (P1)'$ if and only if $E_0(A+K)\subseteq E_0(A)$.

\

\noindent (ii) If $\Phi^{\iso}_{Bw}(A)=\emptyset=\Phi^{\iso}_{Bw}(A+K)$, then (since $\Pi(A)=\Pi_0(A)$, $\Pi^a(A)=\Pi_0(A)\subseteq\Pi_0^a(A)$ and $\Pi^a_0(A)\subseteq E^a_0(A)\subseteq E^a(A)$) 
$$A\in (P2)\Longleftrightarrow E^a(A)=\Pi(A)=\Pi^a(A)\Longrightarrow E^a_0(A)=\Pi_0(A)=\Pi^a_0(A)\Longleftrightarrow A\in (P2)'$$ implies $$E^a_0(A)=\Pi_0(A)=\Pi_0(A+K)\subseteq E^a_0(A+K).$$
If, instead, $\Phi^{\iso}_{uBw}(A)=\emptyset=\Phi^{\iso}_{uBw}(A+K)$, then $\sigma_{Bw}(A)=\sigma_w(A)$ implies 
$\Pi(A)=\Pi_0(A)=\Pi_0(A+K)$, $\Pi^a(A)=\Pi^a_0(A)$ and $\Pi^a_0(A)\subseteq E^a_0(A)\subseteq E^a(A)$. Hence 
$$A\in (P2)\Longleftrightarrow E^a(A)=\Pi(A)=\Pi^a(A)\Longleftrightarrow E^a_0(A)=\Pi_0(A)=\Pi^a_0(A) (\Longleftrightarrow A\in (P2)')$$ implies $$E^a_0(A)=\Pi_0(A+K)=\Pi^a_0(A)=\Pi^a_0(A+K)\subseteq E^a_0(A+K).$$ In either case,  $A\in(P2)\Longrightarrow A\in (P2)'$  and $A\in (P2)\Longrightarrow A+K\in (P2)'$ if and only if $E^a_0(A+K)\subseteq E^a_0(A)$.\end{demo}

The hypotheses of the theorem are not sufficient to guarantee $E(A+K)=E_0(A+K)$, or, $E^a(A+K)=E^a_0(A+K)$ (see Example \ref{ex4}). A sufficient condition ensuring $E_0(A+K)\subseteq E_0(A)$ in (i) and (ii)above is that the operator $A$ is finitely a-polaroid. This follows since $\lambda\in E_0(A+K)$ or $\lambda\in E^a(A+K)$ implies $\lambda\in\iso\sigma_a(A)$, and the hypothesis $A$ is finitely 
a-polaroid implies $\lambda\in\Pi_0(A)$ ($= E_0(A)$ in case (i) and $= E^a_0(A)$ in case (ii)). Observe that the hypothesis $\iso\sigma_a(A)\cap\sigma_w(A)=\emptyset$ guarantees both $\Phi^{\iso}_{Bw}(A)=\emptyset$ and $A$ ia a-polaroid, and the hypothesis $\iso\sigma_{aw}(A)=\emptyset$ guarantees both $\Phi^{\iso}_{uBw}(A)=\emptyset$ and $A$ ia left polaroid. 

\vskip4pt 

Hypotheses $\iso\sigma_a(A)=\iso\sigma_a(A+K)$ and $\Phi^{\iso}_{uBw}(A)=\Phi^{\iso}_{uBw}(A+K)$, where $A, K\in\b$ and $K$ is compact, imply $$\sigma_{Bw}(A)\setminus\sigma_{uBw}(A)=\sigma_w(A)\setminus\sigma_{aw}(A)= \sigma_w(A+K)\setminus\sigma_{aw}(A+K)
=\sigma_{Bw}(A+K)\setminus\sigma_{uBw}(A+K)$$ and $$\iso\sigma_a(A)\cap\{\sigma_{Bw}(A+K)\setminus\sigma_{Bw}(A)\}=\iso\sigma_a(A)\cap\{\sigma_{w}(A+K)\setminus\sigma_{w}(A)\}=\emptyset.$$ Observe here that if $A\in (P2)$ and $E^a(A+K)\subseteq E^a(A)$, then $E^a(A)=\iso\sigma_a(A)\cap\sigma_{uBw}(A)^{\c}$ and hence $\lambda\in E^a(A+K)$ implies $\lambda\notin\sigma_{uBw}(A)$. The following theorem is an analogue of Theorem \ref{thm2} for operators $A\in (P1)$ or  $(P2)$ such that $A+K\in (P1)$ or (respectively) $(P2)$.
\begin{thm}
	\label{thm3} Given operators $A, K\in\b$ with $K$ compact, if $\iso\sigma_a(A)=\iso\sigma_a(A+K)$, $\sigma_{Bw}(A)\cap\sigma_{uBw}(A)^{\c}=\sigma_{Bw}(A+K)\cap\sigma_{uBw}(A+K)^{\c}$ and $\iso\sigma_a(A)\cap\{\sigma_{Bw}(A+K)\setminus\sigma_{Bw}(A)\}=\emptyset$, then a sufficient condition for $$ A\in (Pi)\Longrightarrow A+K\in (Pi), \hspace{4mm} i=1,2,$$ is that $\iso\sigma_a(A)\cap\sigma_{uBw}(A)=\emptyset$.
\end{thm}
\begin{demo} Since $A\in (P1)$ if and only if $E(A)=\Pi^a(A)=\Pi(A)$, and $A\in (P2)$ if and only if $E^a(A)=\Pi(A)=\Pi^a(A)$, the hypothesis $A\in (Pi)$, $i=1,2$, implies $\Pi^a(A)=\Pi(A)$. Hence, if $\iso\sigma_a(A)=\iso\sigma_a(A+K)$ and $\sigma_{Bw}(A)\cap\sigma_{uBw}(A)^{\c}=\sigma_{Bw}(A+K)\cap\sigma_{uBw}(A+K)^{\c}$, then
	\begin{eqnarray*}
	\emptyset &=& \Pi^a(A)\setminus\Pi(A)=\{iso\sigma_a(A)\cap\sigma_{uBw}(A)^{\c}\}\cap\{iso\sigma_a(A)\cap\sigma_{Bw}(A)^{\c}\}^{\c}\\ &=&
	 \{\iso\sigma_a(A)\cap\sigma_{uBw}(A)^{\c}\cap\iso\sigma_a(A)^{\c}\}\cup\{\iso\sigma_a(A)\cap\sigma_{uBw}(A)^{\c}\cap\sigma_{Bw}(A)\}\\ &=& \iso\sigma_a(A)\cap\{\sigma_{Bw}(A)\setminus\sigma_{uBw}(A)\}\\ &=& \iso\sigma_a(A+K)\cap\{\sigma_{Bw}(A+K)\setminus\sigma_{uBw}(A+K)\}\\ &=& \{\iso\sigma_a(A+K)\cap\sigma_{uBw}(A+K)^{\c}\}\cap\{\iso\sigma_a(A+K)\cap\sigma_{Bw}(A+K)^{\c}\}^{\c}\\ &=& \Pi^a(A+K)\setminus\Pi(A+K),
	\end{eqnarray*} i.e., $\Pi^a(A+K)\subseteq\Pi(A+K)$. Since $\Pi(A+K)\subseteq\Pi^a(A+K)$ always, $$\Pi(A+K)=\Pi^a(A+K).$$ Again, since \begin{eqnarray*} \Pi(A)\setminus\Pi(A+K) &=& \Pi(A)\cap\{\iso\sigma_a(A+K)\cap\sigma_{Bw}(A+K)^{\c}\}^{\c}\\ &=& \{\iso\sigma_a(A)\cap\sigma_{Bw}(A)^{\c}\}\cap\{\iso\sigma_a(A)\cap\sigma_{Bw}(A+K)^{\c}\}^{\c}\\ &=&\iso\sigma_a(A)\cap\{\sigma_{Bw}(A+K)\setminus\sigma_{Bw}(A)\}=\emptyset,	
	\end{eqnarray*} we must have $$\Pi(A)\subseteq \Pi(A+K).$$ Consider now a $\lambda\in E(A+K)$. If $\iso\sigma_a(A)\cap\sigma_{uBw}(A)=\emptyset$, then $A$ is left polaroid, hence $\lambda\in E(A+K)$ implies $$\lambda\in\iso\sigma_a(A)\cap\sigma_{uBw}(A)^{\c}=\Pi^a(A)=\Pi(A)\Longrightarrow E(A+K)\subseteq E(A)$$ and, since $\Pi(A+K)\subseteq E(A+K)\subseteq E(A)=\Pi(A)\subseteq\Pi(A+K)$, $$E(A+K)=\Pi^a(A+K)\Longleftrightarrow A+K\in (P1).$$
	Considering, instead, a $\lambda\in E^a(A+K)$, the above argument implies $$\lambda\in\Pi^a(A)=\Pi(A)\subseteq\Pi(A+K)$$ and hence, since $\Pi(A+K)=\Pi^a(A+K)\subseteq E^a(A+K)$, $$E^a(A+K)=\Pi(A+K)\Longleftrightarrow A+K\in (P2).$$ This completes the proof.\end{demo}
If $A\in (P1)$, then the hypotheses of Theorem \ref{thm3} imply $E(A)=\Pi^a(A)=\Pi(A)\subseteq \Pi(A+K)=\Pi^a(A+K)\subseteq E(A+K)$; similarly, if $A\in (P2)$, then the hypotheses of Theorem \ref{thm3} imply $E^a(A)=\Pi(A)=\Pi^a(A)\subseteq \Pi(A+K)=\Pi^a(A+K)\subseteq E^a(A+K)$. Hence a necessary and sufficient condition for $A\in (P1)$ implies $A+K\in (P1)$ (resp. $A\in (P2)$ implies $A+K\in (P2)$) in Theorem \ref{thm3} is that $E(A+K)\subseteq \Pi(A)$ (resp., $E^a(A+K)\subseteq\Pi^a(A)$).
\begin{cor}
	\label{cor1}If $A, K\in\b$ satisfy the hypotheses of Theorem \ref{thm3}, then \begin{eqnarray*}
		& & A\in (P1)\Longrightarrow A+K\in (P1)\Longleftrightarrow E(A+K)\cap\sigma_{Bw}(A)=\emptyset,\hspace{2mm}\mbox{and}\\ & & 
	 A\in (P2)\Longrightarrow A+K\in (P2)\Longleftrightarrow E^a(A+K)\cap\sigma_{uBw}(A)=\emptyset.
	\end{eqnarray*}
\end{cor}

\begin{demo} A straightforward consequence of the facts that $E(A+K)\subseteq\Pi(A)$ if and only if $E(A+K)\cap\sigma_{Bw}(A)=\emptyset$ and  $E^a(A+K)\subseteq\Pi^a(A)$ if and only if $E^a(A+K)\cap\sigma_{uBw}(A)=\emptyset$.
	\end{demo}
We conclude this section with a remark on Hilbert space operators.
\begin{rema}
	\label{rema00}
{\em Given a Hilbert space operator $A\in\B$, there always exists a compact operator $K\in\B$ such that $\sigma_p(A+K)=\Phi_{sf}^+(A)=\{\lambda\in \sigma(A):A-\lambda$ is semi-Fredholm of $\ind(A-\lambda)>0\}=\Phi^+_{sf}(A+K)$ \cite[Proposition 3.4]{HTW}. Consider a $\lambda\in E(A+K)=\sigma_p(A+K)\cap\i\sigma(A+K)$, or, $\lambda\in E^a(A+K)=\sigma_p(A+K)\cap\i\sigma_a(A+K)$. Since $A+K$ has SVEP at $\lambda\in \Phi_{sf}(A+K)$ implies $\ind(A+K-\lambda)\leq 0$ \cite{A}, we have $E(A+K)=E^a(A+K)=\emptyset$. Hence $E(A+K)=\Pi(A+K)=\Pi^a(A+K)=E^a(A+K)=\emptyset$, and $A+K\in (P1)\wedge (P2)$. Conclusion: Given a Hilbert space operator $A\in\B$, there always exists a compact operator $K\in\B$ such that $A+K\in (P1)\wedge (P2)$.}\end{rema} In the absence of similar results for perturbed Banach space operators, a corresponding remark does not seem possible for Banach space operators.

\section {\sfstp Examples: Analytic Toeplitz operators and operators satisfying the abstract shift condition}
If we let $\Omega$ denote the normalized arc length measure on $\partial{\D}$ and let $H^2=H^2(\partial{D})$ denote the Hardy space of analytic square summable (with respect to $\Omega$) functions, then the Toeplitz operator $T_f$ with symbol $f$ is the operator in $B(H^2)$ defined by $$T_f(g)={\cal P}(fg),\hspace{2mm}g\in H^2,$$ where ${\cal P}$ is the orthogonal projection of $L^2(\partial{\D},\Omega)$ onto $H^2$. The operator $T_f$ is analytic Toeplitz if $f\in H^{\infty}(\partial{\D})$. (We assume in the following that $f\neq$ a constant.)
\
If $A\in B(H^2)$ is an analytic Toeplitz operator, then $\sigma(A)=\sigma_w(A)$ is a connected set, $A$ (satisfies Bishop's property $(\beta)$ and so) has SVEP \cite{LN}, and $A$ has no eigenvalues \cite[Page 139]{Hal}. Hence $$E(A)=E^a(A)=\Pi^a(A)=\Pi(A)=\emptyset \Longrightarrow A\in (P1)\wedge (P2).$$ The connected property of $\sigma_w(A)$ implies that $A+K$ is polaroid for all compact operators $K\in B(H^2)$ \cite[Theorem 6.4]{DK}; the connected property of $\sigma_w(A)$ also implies that $$\sigma_w(A)=\sigma_{Bw}(A)=\sigma_{Bw}(A+K)=\sigma_w(A+K)$$ (consequently, $\Pi(A+K)=\Pi_0(A+K)$) for all compact $K\in B(H^2)$.

\

For $A, K\in B(H^2)$, $A$ analytic Toeplitz and $K$ compact, assume that $E(A+K)\neq\emptyset$ and consider a $\lambda\in E(A+K)$. Since $A+K$ is polaroid, $\lambda\in\Pi(A+K)-\Pi_0(A+K)$ and hence (since $\Pi(A+K)\subset E(A+K)$ always) $E(A+K)=E_0(A+K)=\Pi_0(A+K)=\Pi(A+K)$. Recall that $\Pi^a_0(A+K)=\Pi_0(A+K)$ if and only if $\iso\sigma_a(A+K)\cap\{\sigma_w(A+K)\setminus\sigma_{aw}(A+K)\}=\iso\sigma_a(A+K)\cap\{\sigma_w(A)\setminus\sigma_{aw}(A)\}=\emptyset$. Hence, if we now assume that $\iso\sigma_a(A)=\iso\sigma_a(A+K)$, then $$\lambda\in\iso\sigma_a(A+K)\cap\{\sigma_w(A)\setminus\sigma_{aw}(A)\}\Longrightarrow \lambda\in E^a_0(A)\cap\sigma_w(A),$$
a contradiction since $A$ has no eigenvalues. Conclusion: 
\vskip3pt {\em Given operators $A, K\in B(H^2)$, with $A$ analytic Toeplitz and $K$ compact, if $\iso\sigma_a(A)=\iso\sigma_a(A+K)$, then $A+K\in (P1)'$.}
\vskip3pt\noindent We do not know if the hypothesis $\i\sigma_a(A)=\iso\sigma_a(A+K)$ is sufficient to guarantee $A+K\in (P2)'$. A hypothesis guaranteeing $A+K\in (P1)\wedge (P2)$ is given by the following:
\begin{thm}
	\label{thm4} If $A, K\in B(H^2)$, where $A$ is analytic Toeplitz and $K$ is  compact, satisfy $E^a(A+K)\cap\sigma_{w}(A)=\emptyset$, then  $A+K\in (P1)\wedge (P2)$. 
\end{thm} \begin{demo}  It is clear from the above that if $A$ is analytic Toeplitz and $K$ is compact, then $\sigma_{Bw}(A+K)=\sigma_w(A+K)=\sigma_w(A)=\sigma_{Bw}(A)$ and $E(A+K)=\Pi(A+K)=\Pi_0(A+K)=E_0(A+K)$. Since $\Pi_0(A+K)\subseteq\Pi^a_0(A+K)\subseteq E^a_0(A+K)$ and $\Pi_0(A+K)\subseteq E_0(A+K)\subseteq E^a_0(A+K)$, it follows that $E_0(A+K)=\Pi_0(A+K)=\Pi^a_0(A+K)=E^a_0(A+K)=E^a(A+K)$ if and only if $E^a(A+K)\setminus\Pi_0(A+K)=\emptyset$. We have: $$E^a(A+K)\setminus\Pi_0(A+K)=E^a(A+K)\cap\{\iso\sigma_a(A+K)\cap\sigma_w(A+K)^{\c}\}^{\c}=E^a(A+K)\cap\sigma_w(A),$$
which implies $$E^a(A+K)\setminus\Pi_0(A+K)=\emptyset\Longleftrightarrow E^a(A+K)\cap\sigma_w(A)=\emptyset.$$ This completes the proof.
\end{demo} The sufficient condition of the theorem is necessary too: For if $A+K\in (P1)\wedge (P2)$, then $E^a(A+K)=\Pi(A+K)=\Pi_0(A+K)=E^a_0(A+K)$, and hence $E^a(A+K)\setminus\Pi_0(A+K)=E^a(A+K)\cap\{\iso\sigma_a(A+K)\cap\sigma_w(A+K)^{\c}\}^{\c}=E^a(A+K)\cap\sigma_w(A)=\emptyset$.

\vskip6pt
{\bf Operators satisfying the ``abstract shift condition"} $A\in\b$ satisfies the {\em abstract shift condition}, $A\in (ASC)$, if 
$A^{\infty}(\X)=\cap_{n=1}^{\infty}A^n(\X)=\{0\}$ \cite{LN}. Operators $A\in (ASC)$ satisfy the properties that $\sigma(A)$ is 
connected (so that either $\iso\sigma(A)=\emptyset$, or, $\sigma(A)=\{0\}$ in which case $A$ is quasinilpotent), $\alpha(A-\lambda)=0$ for all non-zero $\lambda\in\sigma(A)$ (so that $A$ has SVEP) and $\sigma(A)=\sigma_w(A)$ \cite{{A},{LN}}. 
If we let $r(A)$ denote the spectral radius of $A$ and define $i(A)$ by $$i(A)=\lim_{n\rightarrow\infty}\{\kappa(A^n)\}^{\frac{1}{n}}
=\sup_{n\rightarrow\infty}\{\kappa(A^n)\}^{\frac{1}{n}},$$ where $$\kappa(A)=\inf\{||Ax||:x\in\X, ||x||=1\}$$ denotes the 
lower bound of $A$, then $\D(0,i(A))\subseteq\sigma(A)$. {\em We assume henceforth that $A$ is not quasinilpotent and $i(A)=r(A)$ 
for operators $A\in (ASC)$.} Given a compact operator $K\in\b$, we prove in the following that $A+K\in (P1)\vee (P2)$ (`inclusive' or) if and only if $\iso\sigma_a(A+K)\cap{\eta}'\sigma_{aw}(A)=\emptyset$, where ${\eta}'\sigma_{aw}(A)$ denotes the bounded component of the complement of $\sigma_{aw}(A)$ in $\sigma_w(A)$.

\vskip4pt An important subclass of 	$(ASC)$ operators is that of {\em weighted right shift operators} $A$, $A\in (WRS)$, in $B(\ell^p)$; $\ell^p=\ell^p(\N), 1\leq p<\infty$. It is well known (see \cite{{A},{LN}} and some of the argument above) that 
\begin{eqnarray*} & & \sigma(A)=\sigma_w(A)=\sigma_{Bw}(A)=\overline{{\D}(0,r(A))}, E(A)=E^a(A)=\emptyset,\\ & & \sigma_a(A)=\sigma_{aw}(A)=\sigma_{uBw}(A)={\partial{\D}(0,r(A))} \end{eqnarray*} for  operators $A\in (ASC)$ (recall: $A$ is non-quasinilpotent and $i(A)=r(A)$),  and \begin{eqnarray*} & & \sigma(A)=\sigma_w(A)=\sigma_{Bw}(A)=\overline{{\D}(0,r(A))}, E(A)=E^a(A)=\emptyset,\\ & & \sigma_a(A)=\sigma_{aw}(A)=\sigma_{uBw}(A)=\{\lambda:i(A)\leq \lambda\leq r(A)\}\end{eqnarray*} for operators $A\in (WRS)$. It is clear that {\em $A\in (P1)\wedge (P2)$ for operators $A\in (ASC)\vee (WRS)$}.
\begin{thm}
	\label{thm5}Given an operator $A\in (ASC)\vee(WRS)$, and a compact operator $K$ such that $K\in \b$ if $A\in (ASC)$ and $K\in B(\ell^p)$ if $A\in (WRS)$, $A+K\in (P1)\vee (P2)$, inclusive or, if and only if $\iso\sigma_a(A+K)\cap\{\lambda:0\leq |\lambda|<i(A)\}=\emptyset$. 
\end{thm}
\begin{demo} If $A\in (ASC)\vee(WRS)$, then $\iso\sigma_w(A)=\iso\sigma_{aw}(A)=\emptyset$ implies that $A+K$ is both polaroid and left-polaroid (see Theorem \ref{thm1}). Consequently, $$\lambda\in E(A+K)\Longrightarrow \lambda\in\Pi(A+K),\hspace{2mm} \mbox{hence}\hspace{2mm}E(A+K)=\Pi(A+K)$$ and  $$\lambda\in E^a(A+K)\Longrightarrow \lambda\in\Pi^a(A+K), \hspace{2mm}\mbox{hence}\hspace{2mm}E^a(A+K)=\Pi^a(A+K).$$ Thus \begin{eqnarray*} & & A+K\in (P1)\vee(P2)\Longleftrightarrow\Pi^a(A+K)=\Pi(A+K)\\ &\Longleftrightarrow& \Pi^a(A+K)\subseteq\Pi(A+K)\Longleftrightarrow \Pi^a(A+K)\setminus\Pi(A+K)=\emptyset\\ &\Longleftrightarrow&\i\sigma_a(A+K)\cap\{\sigma_{Bw}(A+K)\setminus\sigma_{uBw}(A+K)\}=\emptyset\\ &\Longleftrightarrow& \iso\sigma_a(A+K)\cap\{\sigma_w(A+K)\cap\sigma_{aw}(A+K)^{\c}\}\\ & & =\iso\sigma_a(A+K)\cap\{\lambda:0\leq|\lambda|<i(A)\}=\emptyset 
	\end{eqnarray*}(where $i(A)=r(A)$ if $A\in (ASC)$.)\end{demo}

\vskip4pt\noindent{\bf Operators $f(A)$.}\hspace{6mm} Let $f\in\h$, where $A\in (ASC)$ or $WRS)$ or $A$ is an analytic Toeplitz operator. (Recall: If $A\in (ASC)$, then $i(A)=r(A)$ and $A$ is not quasinilpotent.) Since $A$ has SVEP (everywhere) and $\sigma_w(A)=\sigma(A)$,
\begin{eqnarray*} & & \sigma(f(A))=f(\sigma(A))=f(\sigma_w(A))=\sigma_w(f(A));\\ & & f(A) \hspace{2mm}\mbox{is polaroid and}\hspace{2mm}E(f(A))=\Pi(f(A)).
\end{eqnarray*}  Recall that $\sigma_p(A)=\emptyset$: We claim that $\sigma_p(f(A))=\emptyset$. For suppose there exists a $\lambda\in\sigma_p(f(A))$. Then there exists a $\mu\in\sigma(A)$ such that $$f(A)-\lambda=f(A)-f(\mu)=(A-\mu)^{\alpha}p(A)g(A)$$
for some integer $\alpha>0$, a polynomial $p(z)$ such that $p(\mu)\neq 0$ and an analytic function $g(z)$ which does not vanish on $\sigma(A)$. But then $(f(A)-\lambda)x=0$, $x\neq 0$, implies $\mu\in\sigma_p(A)$ -- a contradiction. This proves our claim. The fact that $\sigma_p(f(A))=\emptyset$ implies $E^a(f(A))=\emptyset$ ensures (since $\Pi(f(A))\subseteq\Pi^a(f(A))\subseteq E^a(f(A))$) that $$E(f(A))=\Pi(f(A))=\Pi^a(f(A))= E^a(f(A))\Longleftrightarrow f(A)\in (P1)\wedge (P2).$$
Consider now operators $A, K\in B(H^2)$ such that $A$ is analytic Toeplitz and $K$ is compact. Given $f\in\h$, $\i\sigma(f(A))=\i f(\sigma(A))$, $f(A+K)$ is polaroid and hence $$E(f(A+K))=\Pi(f(A+K)).$$  Assume further that $f$ is injective  and $\i\sigma_a(A+K)=\i\sigma_a(A)$. Then $\i\sigma_a(f(A+K))=\i\sigma_a(f(A))$, hence (since $A$ has no eigenvalues) \begin{eqnarray*} & & \i\sigma_a(f(A+K))\cap\{\sigma_w(f(A+K))\setminus\sigma_{aw}(f(A+K))\}\\ &=& f(\i\sigma_a(A+K)\cap\{\sigma_w(A+K)\setminus\sigma_{aw}(A+K)\})\\ &=& f(\i\sigma_a(A)\cap\{\sigma_w(A)\setminus\sigma_{aw}(A)\})=f(\Pi^a_0(A)\cap\sigma_w(A))\\ &=& \emptyset
\end{eqnarray*}  (since $\Pi^a_0(A)=\Pi_0(A)$). Thus: \begin{pro}
\label{pro1} If $f\in\h$ is injective, then $f(A+K)\in (P1)$ for analytic Toeplitz operators $A\in B(H^2)$ perturbed by a compact operator $K\in B(H^2)$ such that $\i\sigma_a(A+K)=\i\sigma_a(A)$.
\end{pro}

The following theorem, an analogue of Theorem \ref{thm4}, gives a necessary and sufficient condition for $f(A+K)\in (P1)\wedge (P2)$.
\begin{thm} 
	\label{thm6} Given operators $A, K\in B(H^2)$, where $A$ is analytic Toeplitz and $K$ is compact, and an injective function $f\in\h$, $f(A+K)\in (P1)\wedge (P2)$ if and only if $E^a(A+K)\cap\sigma_w(A)=\emptyset$.
\end{thm}
\begin{demo} If the operators $A, K$ and the function $f$ are as in the statement of the theorem, then $\i\sigma_x(f(A+K))=\i\sigma_x(f(A))$, $\sigma_x=\sigma$ or $\sigma_a$, $f(A+K)$ is polaroid (hence $E(f(A+K))=\Pi(f(A+K))=f(\Pi(A+K))$) and $f(A+K)$ is left polaroid (so that $\Pi^a(f(A+K))=f(\Pi^a(A+K))$). Consequently, $f(A+K)\in (P1)\wedge (P2)$ if and only if \begin{eqnarray*}  & & E(f(A+K))=\Pi(f(A+K))=\Pi^a(f(A+K))=E^a(f(A+K))\\ &\Longleftrightarrow& E^a(f(A+K))\setminus\Pi(f(A+K))=\Pi^a(f(A+K))\setminus\Pi(f(A+K))=\emptyset.
	\end{eqnarray*} Recalling that $\sigma_{Bw}(A+K)=\sigma_w(A+K)=\sigma_w(A)$, we have 
	\begin{eqnarray*} f(A+K)\in (P1)\wedge (P2) &\Longleftrightarrow& \Pi^a(f(A+K))\setminus\Pi(f(A+K))=\emptyset \\ &\Longleftrightarrow& f(\Pi^a(A+K)\setminus\Pi(A+K))=\emptyset\\ &\Longleftrightarrow& f(\Pi^a(A+K)\cap\sigma_{Bw}(A+K))=\emptyset\\ &\Longleftrightarrow& 
		f(E^a(A+K)\cap\sigma_w(A))=\emptyset\\ &\Longleftrightarrow& E^a(A+K)\cap\sigma_w(A)=\emptyset.\end{eqnarray*}This completes the proof.\end{demo}
	For operators $A\in (ASC)\vee (WRS)$, $\i\sigma_a(A+K)=\i\sigma_a(A)$ ($=\emptyset$) implies $f(A+K)\in (P1)\wedge (P2)$ for all injective $f\in\h$: The hypothesis $\i\sigma_a(A+K)=\emptyset$ may be relaxed.
	\begin{thm}
		\label{thm7} Given operators $A$ and $K$, where $A\in (ASC)\vee (WRS)$ and $K$ is compact, and an injective $f\in\h$, $f(A+K)\in (P1)\vee (P2)$ if and only if $\i\sigma_a(A+K)\cap\{\lambda:0\leq|\lambda|< i(A)\}=\emptyset$.
	\end{thm}
\begin{demo} The injective hypothesis on $f\in\h$ implies $$\i\sigma_x(f(A+K))=f(\i\sigma_x(A+K)), \sigma_x=\sigma\hspace{2mm}\mbox{or}\hspace{2mm}\sigma_a,$$ and $A+K$ (alongwith being polaroid) is left polaroid. Since $$\sigma_{Bw}(f(A+K))=\sigma_w(f(A+K))=f(\sigma_w(A+K))=f(\sigma_{Bw}(A+K))$$ and $$\sigma_{uBw}(f(A+K))=\sigma_{aw}(f(A+K))=f(\sigma_{aw}(A+K))=f(\sigma_{uBw}(A+K)),$$ we have 
	$$E(f(A+K))=\Pi(f(A+K))\hspace{2mm}\mbox{and}\hspace{2mm} E^a(f(A+K))=\Pi^a(f(A+K)).$$ Thus 
	
\begin{eqnarray*}  & & f(A+K)\in (P1)\vee (P2)\Longleftrightarrow \Pi^a(f(A+K))\subseteq\Pi(f(A+K))\\ &\Longleftrightarrow& \i\sigma_a(f(A+K))\cap f\{\sigma_{Bw}(A+K)\setminus\sigma_{uBw}(A+K)\}=\emptyset\\ &\Longleftrightarrow& f(\i\sigma_a(A+K)\cap\{\sigma_w(A+K)\setminus\sigma_{aw}(A+K)\})=\emptyset\\ &\Longleftrightarrow& \i\sigma_a(A+K)\cap\{\sigma_w(A)\cap\sigma_{aw}(A)^{\c}\}=\emptyset\\ &\Longleftrightarrow& \i\sigma_a(A+K)\cap\{\lambda:0\leq|\lambda|< i(A)\}=\emptyset.\end{eqnarray*} (Recall: $i(A)=r(A)$ for $A\in (ASC)$.)\end{demo}

\vskip10pt \noindent\normalsize\rm B.P. Duggal, 8 Redwood Grove,
 London W5 4SZ, United Kingdom.\\
\noindent\normalsize \tt e-mail: bpduggal@yahoo.co.uk

\end{document}